\documentclass[ECP]{ejpecp} 

\SHORTTITLE{Strong and Weighted Matchings}

\TITLE{Strong and Weighted Matchings in Inhomogenous Random Graphs}

\AUTHORS{%
  Ghurumuruhan Ganesan\footnote{Institute of Mathematical Sciences, HBNI, Chennai.
    \EMAIL{gganesan82@gmail.com}}}


\AMSSUBJ{60J10} 

\SUBMITTED{} 
\ACCEPTED{} 

\VOLUME{0}
\YEAR{2020}
\PAPERNUM{0}
\DOI{10.1214/YY-TN}

\ABSTRACT{We equip the edges of a deterministic graph~\(H\) with independent but not necessarily identically distributed weights and study a generalized version of matchings (i.e. a set of vertex disjoint edges) in~\(H\) satisfying the property that end-vertices of any two distinct edges
are at least a minimum distance apart. We call such matchings as strong matchings  and determine bounds on the expectation and variance of the minimum weight of a maximum strong matching. Next, we consider an inhomogenous random graph whose edge probabilities are not necessarily the same and determine bounds on the maximum size of a strong matching in terms of the averaged edge probability. We use local vertex neighbourhoods, the martingale difference method and iterative exploration techniques to obtain our desired estimates.} 

\newcommand{\ind}{1\hspace{-2.3mm}{1}}

\begin{document}



\setcounter{equation}{0}
\renewcommand\theequation{\thesection.\arabic{equation}}
\section{Introduction}
Matchings in graphs is an important object of study from both theoretical and application perspectives.  One of the most well-studied aspect of matchings from the probabilistic perspective is that of minimum weight matchings~\cite{ald} : Given a complete bipartite graph on~\(n+n\) vertices and equipping each edge with an independent exponential weight, the problem is to determine the minimum weight~\(C(n)\) of a perfect matching. It is well known~\cite{wast2,nair} that~\(\mathbb{E}C(n) = \sum_{k=1}^{n} \frac{1}{k^2}\) and later~\cite{wast} obtained estimates for the expected minimum weight of a matching of given size. Recently~\cite{frie} studied minimum weight matchings in \emph{random} graphs and used the Talagrand concentration inequalities to get the corresponding deviation estimates. 

In the first part of our paper, we study minimum weight of a strong matchings where end-vertices of distinct edges are at a given minimum distance apart. We equip the edges of a deterministic graph~\(H\) with independent but not necessarily identically distributed weights and use local neighbourhood estimates to obtain bounds on the expectation and variance of the minimum weight of a maximum strong matching. We also determine sufficient conditions so that the minimum weight grows linearly with the strong matching number. 



The second main result of our paper is regarding the strong matching number of inhomogenous random graphs where the edge probabilities need not be the same. Matchings in homogenous random graphs where edge probabilities are the same, have been studied extensively and bounds for the matching number are known for a wide range of the edge probability~\cite{boll2}. Any lower bound on the matching number for homogenous graphs can be extended to inhomogenous graphs whose edge probabilities are bounded from below, by the following monotonicity property: If~\(H_1 \subseteq H_2\) are two graphs, then a matching in~\(H_1\) is also a matching in~\(H_2.\)

Induced matchings have been studied in~\cite{clark} (see also references therein) under the name of strong matchings where estimates for the \emph{expected} size of a maximum induced matching in homogenous graphs with constant edge probability is obtained. Recently~\cite{cool} used a combination of second moment method along with concentration inequalities to estimate the largest possible size of induced matchings in homogenous graphs and obtained deviation bounds for a wide range of edge probabilities. A main bottleneck in directly extending the above results to inhomogenous graphs is that induced matchings do not satisfy the monotonicity property enjoyed by the ordinary matchings described in the previous paragraph.  In this paper, we study generalized strong matchings in inhomogenous random graphs and use local neighbourhood bounds described earlier, to estimate the strong matching number in terms of the \emph{averaged} edge probabilities.

The paper is organized as follows: In Section~\ref{weight_strong}, we define strong matchings in graphs and state and prove our first main result, Theorem~\ref{thm_gen}, regarding the minimum weight of a strong matching in a graph equipped with inhomogenous weights. Next in Section~\ref{inhom_strong},  we state and prove Theorem~\ref{gnp_str} regarding the maximum size of a strong matching in inhomogenous random graphs.

\setcounter{equation}{0}
\renewcommand\theequation{\thesection.\arabic{equation}}
\section{Weighted Strong Matchings}\label{weight_strong}
Let~\(H = (V,E)\) be any graph containing at least one edge. A path~\(\pi\) in~\(H\) is a sequence of edges~\((e_1,\ldots,e_l)\) such that~\(e_{i} = (a_{i},b_i)\) and~\(e_{i+1} = (a_{i+1},b_{i+1})\) share a common endvertex~\(b_i = a_{i+1}\) for~\(1 \leq i \leq l-1.\) The length of~\(\pi\) is equal to~\(l,\) the number of edges in~\(\pi\) and the vertices~\(a_1\) and~\(b_l\) are said to be connected by the path~\(\pi.\)  The distance between two vertices~\(a\) and~\(b\) is defined to be the minimum length of a path connecting~\(a\) and~\(b.\)

A set of vertex disjoint edges~\({\cal W} = \{e_1,\ldots,e_l\}\) in~\(H\) is said to be a \emph{matching} of size~\(l.\)
\begin{definition}\label{def_one} Let~\({\cal W} =\{e_1,\ldots,e_l\}\) be a matching of a graph~\(H.\) For integer~\(k \geq 0,\) we say that~\({\cal W}\) is a~\(k-\)\emph{strong matching} if the following property holds for any pair of edges~\(e_i \neq e_j:\) There does not exist a path~\(\pi\) in~\(H\) containing~\(l \leq k\) edges that connects an endvertex of~\(e_i\) with an endvertex of~\(e_j.\)

The~\(k-\)\emph{strong matching number}~\(\nu_k(H)\) of the graph~\(H\) is the size of a maximum~\(k-\)strong matching in~\(H.\)
\end{definition}
For~\(k=0\) the above definition reduces to the usual definition of matching as described prior to Definition~\ref{def_one} and for~\(k=1\) this coincides with the concept of induced/strong matchings described in the introduction. We retain the terminology strong matchings and introduce the term~\(k\) as a further generalization.

We now equip each edge of~\(H\) with random weights and estimate the minimum weight of a strong matching. Let~\(w(e) \geq  0\) be the  random weight assigned to the edge~\(e \in H.\) The edge weights are independent but not necessarily identically distributed and we define~\(M_k(H)\) to be the minimum weight of a maximum~\(k-\)strong matching of~\(H.\) We have the following properties regarding the mean and variance of the minimum weight~\(M_k(H).\)
We say that an edge~\(e =(u,v)\) of a graph~\(H\) is isolated if no other edge in~\(H\) shares an endvertex with~\(H.\)
\begin{theorem}\label{thm_gen} Suppose the graph~\(H\) has no isolated edges.\\
\((a)\) Letting
\begin{equation}\label{mu_def}
\mu := \max_{e \in H} \mathbb{E}w(e) \text{ and } \mu_2 := \max_{e \in H} \mathbb{E}w^2(e),
\end{equation}
we have that
\begin{equation}\label{mk_mean2}
\mathbb{E}M_k(H) \leq \mu \nu_k(H) \text{ and } var(M_k(H)) \leq 4\mu_2\nu_k(H),
\end{equation}
\((b)\) Suppose the edge weights are stochastically dominated by distributions~\(F_1\) and~\(F_2\) in the sense that
\begin{equation}\label{sto_dom}
0 < F_1(x) \leq \mathbb{P}(w(e) \leq x) \leq F_2(x) \text{ for each } x > 0 \text{ and each edge }e
\end{equation}
and~\(F_2(x) \longrightarrow 0\) as~\(x \rightarrow 0.\) If~\(\Delta\) is the maximum degree of a vertex in~\(H,\) then
\begin{equation}\label{mk_mean}
\mathbb{E}M_k(H) \geq C \nu_k(H)
\end{equation}
where~\(C > 0\) is  a positive constant that depends only on~\(\Delta,k,F_1\) and~\(F_2.\)\\
\end{theorem}
We have the following remarks:\\
\underline{\emph{Remark 1}}: The condition~(\ref{sto_dom}) in~\((b)\) is satisfied if for example, the weights are exponentially distributed with unit mean. The estimate~(\ref{mk_mean}) then implies that if~\(H\) is a \emph{bounded} degree graph, then the expected minimum weight of a strong matching grows linearly with the strong matching number. Moreover, from the variance estimate in~(\ref{mk_mean2}) we see that the minimum weight of a strong matching is concentrated around its mean.\\
\underline{\emph{Remark 2}}: The bounded degree condition is important since in the proof of~(\ref{mk_mean}) below, we see that the constant~\(C\) is such that~\(C \longrightarrow 0\) as the maximum vertex degree~\(\Delta \longrightarrow \infty.\)

\subsection*{Proof of Theorem~\ref{thm_gen}}
\emph{Proof of Theorem~\ref{thm_gen}~\((a)\)}: Let~\({\cal W}\) be any deterministic maximum~\(k-\)strong matching of~\(H\) containing~\(\nu_k(H)\) edges and let~\(M_k({\cal W}) = \sum_{e \in {\cal W}} w(e)\) be the weight of~\({\cal W}.\) We then have that~\(\mathbb{E}M_k(H) \leq \sum_{e \in {\cal W}} \mathbb{E}w(e) \leq \mu \cdot \nu_k(H)\) where~\(\mu\) is as in the statement of the Theorem. This proves the first bound in~(\ref{mk_mean2}).

For the variance bound, we use the martingale difference method analogous to~\cite{kest}. Let~\(e_1,\ldots,e_q\) be the set of edges of the graph~\(H\) and for~\(1 \leq j \leq q,\) let~\({\cal F}_j = \sigma\left(\{w(e_k)\}_{1 \leq k \leq j}\right)\) denote the sigma field
generated by the weights of the edges~\(\{e_k\}_{1 \leq k \leq j}.\) We define the martingale difference
\begin{equation}\nonumber
\zeta_j = \mathbb{E}(M_k(H) \mid {\cal F}_j) - \mathbb{E}(M_k(H) \mid {\cal F}_{j-1}),
\end{equation}
and get that~\(M_k(H) -\mathbb{E}M_k(H) = \sum_{j=1}^{q} \zeta_j.\) By the martingale property we then have
\begin{equation} \label{var_exp}
var(M_k(H))  = \mathbb{E}\left(\sum_{j=1}^{q} \zeta_j\right)^2 = \sum_{j=1}^{q} \mathbb{E}\zeta_j^2.
\end{equation}

To evaluate~\(\mathbb{E}\zeta_j^2\) we rewrite~\(\zeta_j = \mathbb{E}(Q(e_j) \mid {\cal F}_j),\)
where~\(Q(e_j) := M_{k}(w(e_j)) - M_{k}(t(e_j))\) has the following terminology: the random variable~\(t(e_j)\) is an independent copy of~\(w(e_j)\) which is also independent of~\(\{w(e_l)\}_{1 \leq l \neq j \leq q}\) and~\(M_{k}(w(e_j))\) and~\(M_{k}(t(e_j))\) are the weights of the minimum weight maximum~\(k-\)strong matchings~\({\cal W}\) and~\({\cal T}\) obtained respectively with edge weights~\(\{w(e_k)\}_{1 \leq k \leq q}\) and~\(\{w(e_k)\}_{1 \leq k \neq j \leq q} \cup \{t(e_j)\}.\)

From the above paragraph we get that
\begin{equation}\label{hj2}
\zeta_j^2 \leq \left(\mathbb{E}(Q(e_j)\mid{\cal F}_j)\right)^2 \leq \mathbb{E}(\Delta^2(e_j) \mid {\cal F}_j)
\end{equation}
and in what follows we estimate~\(|Q(e_j)|.\) If the weight of the edge~\(e_j\) is decreased from~\(w(e_j)\) to~\(t(e_j) < w(e_j),\) then the corresponding minimum weight~\(M_k(t(e_j)) \leq M_k(w(e_j))\) and the difference~\(M_k(w(e_j)) - M_k(t(e_j)) \leq w(e_j) -t(e_j) \leq w(e_j).\) Also we have that~\(M_k(w(e_j))-M_k(t(e_j))\) is non-zero if and only if~\(e_j \in {\cal T}\) or~\(e_j \in {\cal W}.\) But because~\(t(e_j) < w(e_j)\) we have that~\(\{e_j \in {\cal W}\} \subseteq \{e_j \in {\cal T}\}.\)  Summarizing we have~\[|Q(e_j)| \ind(t(e_j) < w(e_j)) \leq w(e_j) \ind(e_j \in {\cal T})\] and arguing similarly for the case~\(t(e_j) > w(e_j),\) we get that
\begin{equation}\label{del_est_fg}
|Q(e_j)| \leq w(e_j) \ind(e_j \in {\cal T}) + t(e_j) \ind(e_j \in {\cal W}).
\end{equation}

Using~\((a+b)^2 \leq 2(a^2+b^2)\) we get from~(\ref{del_est_fg}) that
\[|Q(e_j)|^2 \leq 2(w^2(e_j) \ind(e_j \in {\cal T}) + t^2(e_j) \ind(e_j \in {\cal W}))\] and taking conditional expectations with respect to the sigma field~\({\cal F}_j\) we get
\begin{equation}\label{del_y}
\mathbb{E}\left(|Q(e_j)|^2 \mid {\cal F}_j\right) \leq 2(I_1+I_2)
\end{equation}
where
\begin{equation}\label{i1_est_one}
I_1 := \mathbb{E}(w^2(e_j) \ind(e_j \in {\cal T}) \mid {\cal F}_j) = w^2(e_j) \mathbb{P}\left(e_j \in {\cal T} \mid{\cal F}_{j-1}\right)
\end{equation}
and
\begin{equation}\label{i2_est_two}
I_2 := \mathbb{E}(t^2(e_j) \ind(e_j \in {\cal W}) \mid {\cal F}_j) = \mathbb{E}t^2(e_j)\mathbb{P}\left(e_j \in {\cal W} \mid {\cal F}_j\right) .
\end{equation}
The final equality in~(\ref{i1_est_one}) is true since the random variable~\(w(e_j)\) is independent of the event~\(e_j \in {\cal T},\) which is determined by the edge weights~\(\{w(e_k)\}_{k \neq j} \cup t(e_j).\) Similarly the final relation in~(\ref{i2_est_two}) follows from the fact that the event~\(\{e_j \in {\cal W}\} \in {\cal F}_j\) is independent of~\(t(e_j).\)

Taking conditional expectations with respect to~\({\cal F}_{j-1}\) in~(\ref{i1_est_one}) and~(\ref{i2_est_two}) we get that
\[\mathbb{E}(I_1 \mid {\cal F}_{j-1}) = \mathbb{E}(I_2 \mid {\cal F}_{j-1}) = \mathbb{E}w^2(e_j) \mathbb{P}\left(e_j \in {\cal T}\mid{\cal F}_{j-1}\right)\]
and so from~(\ref{del_y}) we have~\(\mathbb{E}\left(|Q(e_j)|^2 \mid {\cal F}_{j-1}\right) \leq 4\mathbb{E}w^2(e_j) \mathbb{P}\left(e_j \in {\cal T}\mid{\cal F}_{j-1}\right).\) From~(\ref{hj2}) we therefore have~\(\mathbb{E}(\zeta_j^2\mid{\cal F}_{j-1}) \leq \mathbb{E}\left(|Q(e_j)|^2 \mid {\cal F}_{j-1}\right) \leq 4\mathbb{E}w^2(e_j) \mathbb{P}\left(e_j \in {\cal T} \mid{\cal F}_{j-1}\right)\) and taking expectations we get~\(\mathbb{E}(\zeta_j^2) \leq 4 \mathbb{E}w^2(e_j) \mathbb{P}\left(e_j \in {\cal T}\right) \leq 4\mu_2\mathbb{P}\left(e_j \in {\cal T}\right)\) where~\(\mu_2\) is as in the statement of the Theorem. Summing over~\(j\) and using~(\ref{var_exp}) we then get
\begin{equation}\label{mk_h_lat}
var(M_k(H)) \leq 4\mu_2\sum_{j} \mathbb{P}\left(e_j \in {\cal T}\right) = 4\mu_2\nu_k(H),
\end{equation}
since the expected number of edges in any maximum~\(k-\)strong matching is~\(\nu_k(H).\)~\(\qed\)

\emph{Proof of Theorem~\ref{thm_gen}~\((b)\)}: We first show that there are positive constants~\(\gamma_1\) and~\(\gamma_2\) depending only on~\(\Delta,k,\mu,F_1\) and~\(F_2\) such that
\begin{equation}\label{dev_low}
\mathbb{P}\left(M_k(H) \geq \gamma_1 \cdot \nu_k(H)\right) \geq 1-e^{-\gamma_2 m},
\end{equation}
where~\(m = m(H)\) is the number of edges in~\(H.\) From~(\ref{dev_low}) we get that
\[\mathbb{E}M_k(H) \geq \gamma_1 \cdot \nu_k(H) \cdot (1-e^{-\gamma_2 m}) \geq \gamma_1 \cdot (1-e^{-\gamma_2}) \cdot \nu_k(H)\]
since~\(m \geq 1\) and this obtains the desired lower bound in~(\ref{mk_mean}).

To prove~(\ref{dev_low}), we use a combinatorial result (Lemma~\ref{prop_one} in Appendix) that obtains bounds for the strong matching number in terms of size of local neighbourhoods. Say that~\(e \in H\) is a \emph{bad} edge if its weight~\(w(e) \leq \gamma\) for some constant~\(\gamma > 0\) to be determined later. Letting~\(T_L = \sum_{e} \ind(w(e)  \leq \gamma)\) be the total number of bad edges, we have that\\\(m \cdot \beta_{low} \leq \theta_L = \mathbb{E}T_L \leq m\cdot \beta_{up},\) where
\begin{equation}\label{beta_low}
\beta_{low} := \min_{e} \mathbb{P}(w(e) \leq \gamma) \geq F_1(\gamma)\;\;\text{ and }\;\;\beta_{up}  := \max_{e}\mathbb{P}(w(e) \leq \gamma) \leq F_2(\gamma)
\end{equation}
by~(\ref{sto_dom}). Consequently, using the deviation estimate~(\ref{conc_est_f}) of Lemma~\ref{lem_one} in Appendix with~\(\epsilon = \frac{1}{2},\) we get that~\(\mathbb{P}\left(T_L \geq \frac{3}{2} \cdot m \cdot \beta_{up}\right) \leq 2\exp\left(-\frac{m \cdot \beta_{low}}{16}\right)\) and therefore that
\begin{equation}\label{tl_eq}
\mathbb{P}\left(T_L \geq \frac{3}{2}m\cdot F_2(\gamma)\right) \leq 2\exp\left(-\frac{m \cdot F_1(\gamma)}{16}\right).
\end{equation}

Suppose now that the event~\(T_L \leq \frac{3}{2} m \cdot F_2(\gamma)\) occurs and let~\({\cal U}\) be a maximum~\(k-\)strong matching of minimum weight. There are at most~\(T_L\) bad edges in~\({\cal U}\) and so
\begin{equation}\label{mk_w}
M_k({\cal U}) \geq \left(\nu_k(H) - T_L\right) \cdot \gamma \geq \left(\nu_k(H) - \frac{3}{2} m F_2(\gamma)\right) \gamma \geq \nu_k(H) \cdot \gamma \cdot \left(1-12\cdot \Delta^{k+1} \cdot F_2(\gamma)\right),
\end{equation}
using the maximum degree bound for~\(\nu_k(H)\) (see~(\ref{nu_g_low}) of Lemma~\ref{prop_one} in Appendix). Since~\(F_2(\gamma) \longrightarrow 0\) as~\(\gamma \rightarrow 0\) (see statement of Theorem~\ref{thm_gen}), we choose~\(\gamma > 0\) small so that~\(F_2(\gamma) \leq \frac{1}{24\cdot \Delta^{k+1}}.\) For such a~\(\gamma\) we get from~(\ref{mk_w}) that~\(M_k({\cal U}) \geq \frac{\nu_k(H) \cdot \gamma}{2}\) and so we get from~(\ref{tl_eq}) that~\(\mathbb{P}\left(M_k(H) \geq \frac{\nu_k(H) \cdot \gamma}{2}\right)\geq 1- 2\exp\left(-\frac{m\cdot F_1(\gamma)}{16}\right).\) This proves~(\ref{dev_low}).~\(\qed\)

\setcounter{equation}{0}
\renewcommand\theequation{\thesection.\arabic{equation}}
\section{Inhomogenous Random Graphs}\label{inhom_strong}
In this section, we study strong matching numbers of \emph{random} graphs obtained as follows. Let~\(K_n\) be the complete graph on~\(n\) vertices and let~\(\{X_e\}_{e \in K_{n}}\) be independent random variables indexed by the edge set of~\(K_{n}\) and having distribution
\begin{equation}\label{xe_def}
\mathbb{P}(X_e = 1) = p(e) = 1-\mathbb{P}(X_e = 0)
\end{equation}
for the edge~\(e.\) Let~\(G\) be the random graph formed by the set of all edges~\(e\) satisfying~\(X_e~=~1.\) Because the edge probabilities~\(p(e)\) need not all be the same, we define~\(G\) to be an \emph{inhomogenous} random graph. If~\(p(e) = p\) for all~\(e,\) then~\(G\) is said to be a \emph{homogenous} random graph with edge probability~\(p.\)

We have the following result regarding the strong matching number of inhomogenous random graphs. Throughout constants do not depend on~\(n.\)
\begin{theorem}\label{gnp_str} Let~\(k \geq 3\) an integer and suppose the edge probabilities are of the form~\(p(e) = \frac{h(e)}{n^{\beta}}\) where~\(0 < \beta < 1\) is a constant and~\(\{h(e)\}_{e \in K_n}\) are ``weights" satisfying
\begin{equation}\label{h_cond}
\gamma_1 \cdot n^{\delta_{low}} \leq h(e) \leq \gamma_2 \cdot n^{\delta_{up}} \text{ and } {{n \choose 2}}^{-1}\sum_{e \in K_n} h^{k+2}(e) \leq \gamma_0 \cdot n^{(k+2)\delta_{av}}
\end{equation}
for some constants~\(0 \leq \delta_{low} \leq \delta_{av} \leq \delta_{up} <  \beta\) and~\(\gamma_i,i=0,1,2.\) Suppose
\begin{equation}\label{thet_def}
\theta_{low} := 1-k(1-\beta+\delta_{av}) - 2(\delta_{av}-\delta_{low}) > 0\;\;\text{ and }\;\;k_1(1-\beta+\delta_{up}) < 1
\end{equation}
strictly where~\(k_1 := \frac{k-1}{2}\) if~\(k\) is odd and~\(k_1 := \frac{k-2}{2}\) otherwise. There are positive constants~\(D_i,1 \leq i \leq  3\) such that
\begin{equation}\label{bound_two}
\mathbb{P}\left( D_1 \cdot n^{\theta_{low}} \leq \nu_k(G) \leq D_2\cdot n^{\theta_{up}}\right) \geq 1 - \frac{D_3}{n^{\theta_{low}}}
\end{equation}
for all~\(n\) large, where~\(\theta_{up} := 1-k_1(1-\beta+\delta_{low}).\)
\end{theorem}
In words, we get that with high probability, i.e. with probability converging to one as~\(n~\rightarrow~\infty,\) the maximum size of a~\(k-\)strong matching grows as a power of~\(n.\) Moreover, our bounds for the strong matching are in terms of the \emph{averaged} edge probability parameter~\(\delta_{av}\) and~\(\delta_{low}.\) In particular if~\(\delta_{low}  =\delta_{av} = \delta,\) then~\(\frac{D_1 \cdot n}{n^{k(1-\beta+\delta)}} \leq \nu_k(G) \leq \frac{D_2 \cdot n}{n^{k_1(1-\beta+\delta)}}\) with high probability, where we recall that~\(k_1\) is roughly equal to~\(\frac{k}{2}.\) Extrapolating the results of~\cite{cool} for induced matchings where~\(k=1,\) we    conjecture that~\(\nu_k(G)\) is in fact of the order of~\(\frac{n (\log{n})^{a}}{n^{k(1-\beta+\delta)}}\) with high probability, for some constant~\(a  >0.\)


The proof outline for Theorem~\ref{gnp_str} is as follows. We use local vertex neighbourhood bounds for the strong matching number obtained in Lemma~\ref{prop_one} of Appendix together with variance estimates for weighted random graphs described in Theorem~\ref{thm_gen}, to obtain our lower bound in Theorem~\ref{gnp_str}. For the upper bound on the strong matching number in Theorem~\ref{gnp_str}, we use iterative exploration techniques analogous to~\cite{boll} to first estimate the size of local neighbourhoods of vertices. We then employ an upper bound for the strong matching number again based on local neighbourhoods, obtained in Lemma~\ref{prop_one} of Appendix, to complete the proof of Theorem~\ref{gnp_str}.

\subsection*{Proof of Theorem~\ref{gnp_str}}
\emph{Proof of the lower bound in~(\ref{bound_two})}: We use the estimate~(\ref{nu_g_low}) of Lemma~\ref{prop_one} which states that
\begin{equation}\label{nu_low_k}
\nu_k(G) \geq \frac{m^2}{4\sum_{u \in V} d_1(u)d_{k+1}(u)}
\end{equation}
provided~\(G\) contains no isolated edge. Here~\(d_j(u)\) is the number of vertices at a distance at most~\(j\) from~\(u.\) The expected number of neighbours of any vertex in~\(G\) is at least~\(\frac{\gamma_1 (n-1)}{n^{\beta-\delta_{low}}}\) and so from the deviation estimate~(\ref{conc_est_f}) of Lemma~\ref{lem_one} in Appendix,  we see that each vertex has degree at least~\(\frac{\gamma_1}{2} \cdot n^{1-\beta+\delta_{low}}\) with probability at least\\\(1-e^{-2Cn^{1-\beta+\delta_{low}}}\) for some constant~\(C>0.\) Therefore if~\(E_{iso}\) denotes the event that~\(G\) contains no isolated edge, then by the union bound we have that
\begin{equation}\label{e_iso}
\mathbb{P}(E_{iso}) \geq 1-n \cdot e^{-2Cn^{1-\beta+\delta_{low}}} \geq 1- e^{-Cn^{1-\beta+\delta_{low}}}
\end{equation}
for all~\(n\) large.

To lower bound the number of edges~\(m,\) we use the lower bound for~\(h(e)\) in~(\ref{h_cond}) to get that the expected number of edges in the graph~\(G\) is at least~\({n \choose 2} \cdot \frac{\gamma_1}{n^{\beta-\delta_{low}}} \geq \gamma_1 \frac{n^{2-\beta+\delta_{low}}}{4}\) for all~\(n\) large. Setting~\(\epsilon = \frac{1}{2}\) in the deviation estimate~(\ref{conc_est_f}), we get that the number of edges~\(m\) in~\(G\) satisfies
\begin{equation}\label{m_bound}
\mathbb{P}\left(m \geq \gamma_1 \frac{n^{2-\beta+\delta_{low}}}{8}\right) \geq 1-e^{-Cn^{2-\beta+\delta_{low}}},
\end{equation}
for some constant~\(C > 0.\)

In what follows we find an upper bound for the denominator term~\(\sum_{u \in V} d_1(u)d_{k+1}(u)\) in~(\ref{nu_low_k}). Letting~\(u \sim v\) denote that vertex~\(u\) is adjacent to vertex~\(v,\) we get that the degree of~\(u\) equals~\(d_1(u) = \sum_{v \neq u} \ind(u \sim v),\) where~\(\ind(.)\) is the indicator function.  To compute~\(d_{k+1}(u),\) we use the fact that if a vertex~\(z\) is at a distance~\(l\) from~\(u,\) then there is a path of length~\(l\) containing~\(u\) as an endvertex.
Therefore
\[d_{k+1}(u) \leq \sum_{1 \leq j \leq k+1}\sum_{w_1\neq w_2 \neq \ldots \neq w_{j} \neq u} \ind(u \sim w_1) \prod_{i=1}^{j-1} \ind(w_i \sim w_{i+1})\]
and consequently,
\begin{eqnarray}\label{jk_def}
J_k(n) &:=& \mathbb{E}\sum_{u} d_1(u) d_{k+1}(u) \nonumber\\
&\leq& \sum_{u} \sum_{v \neq u } \sum_{1 \leq j \leq k+1} \sum_{w_1 \neq \ldots \neq w_{j} \neq u} \mathbb{E}f(v,u,w_1,\ldots,w_{j}) \nonumber\\
&=& \sum_{1 \leq j \leq k+1} \sum_{u} \sum_{v \neq u } \sum_{w_1 \neq \ldots \neq w_{j} \neq u} \mathbb{E}f(v,u,w_1,\ldots,w_{j}).
\end{eqnarray}
where
\[f_j  = f(v,u,w_1,\ldots,w_{j}) := \ind(v \sim u) \ind(u \sim w_1)\prod_{i=1}^{j-1} \ind(w_i \sim w_{i+1}).\]

Depending on whether~\(v =w_1\) or not, we split~\(J_k(n)\leq U_k(n) + V_k(n)\) where
\begin{equation}\label{uk_def}
U_k(n) := \sum_{1 \leq j \leq k+1} \sum_{v \neq u \neq w_1}\sum_{u \neq w_1 \neq \ldots \neq w_{j}} \mathbb{E}f(v,u,w_1,\ldots,w_{j})
\end{equation}
and
\begin{equation}\label{vk_def}
V_k(n) := \sum_{1 \leq j \leq k+1} \sum_{u \neq v = w_1 \neq \ldots \neq w_{j}} \mathbb{E}f(v,u,w_1,\ldots,w_{j})
\end{equation}
respectively.

Define a walk to be a sequence of edges~\((e_1,\ldots,e_t)\) where~\(e_i\) and~\(e_{i+1}\) share a common endvertex for~\(1 \leq i \leq t-1.\) The term~\(f_j = f(v,u,w_1,\ldots,w_{j})\) is a product of exactly~\(j+1\) distinct terms if~\(v \neq w_1,\) since in this case the walk~\({\cal W} = (v,u,w_1,\ldots,w_j)\) contains the path~\((u,w_1,\ldots,w_j)\) formed by~\(j\) edges and the additional edge~\((v,u).\) The total number of vertices in the walk~\({\cal W}\) is either~\(j+2\) or~\(j+1\) depending on whether\\\(v \in \{w_2,\ldots,w_j\}\) or not. In any case, the expectation
\begin{equation}\label{efj}
\mathbb{E}f_j = \frac{\prod_{e \in {\cal W}} h(e)}{n^{(j+1)\beta}} \leq \frac{\sum_{e \in {\cal W}} h^{j+1}(e)}{j+1} \cdot \frac{1}{n^{(j+1)\beta}} \leq \sum_{e \in {\cal W}} h^{j+1}(e) \cdot \frac{1}{n^{(j+1)\beta}}
\end{equation}
where the first inequality in~(\ref{efj}) is true by the AM-GM inequality. Substituting~(\ref{efj}) into~(\ref{jk_def}) we then get that
\begin{equation}\label{uk_est2}
U_k(n) \leq \sum_{1 \leq j \leq k+1} \frac{1}{n^{(j+1)\beta}} \sum_{{\cal W}} \sum_{e \in {\cal W}}h^{j+1}(e) .
\end{equation}

The number of walks containing any edge~\((u,v)\) and having~\(l \leq j+2\) vertices is at most~\(C_1 \cdot n^{j}\) for some constant~\(C_1 > 0\) and so each edge of the complete graph~\(K_n\) is counted at most~\(C_1 \cdot n^{j}\) times in the inner double summation in~(\ref{uk_est2}). Thus
\begin{equation}\label{uk_est3}
U_k(n) \leq C_1 \sum_{1 \leq j \leq k+1} \sum_{f \in K_n} h^{j+1}(f) \cdot n^{j-(j+1)\beta}.
\end{equation}
and because~\(j+1 \leq k+2\) and the edge probability ``weights"~\(h(e) \geq \gamma_1\) (see~(\ref{h_cond})), we also have that
\[\sum_{f \in K_n} h^{j+1}(f) \leq C_2 \cdot \sum_{f \in K_n} h^{k+2}(f) \leq C_3 \cdot n^{(k+2)\delta_{av}} n^2\] for some positive constants~\(C_2,C_3,\) by the condition~(\ref{h_cond}). Therefore
\begin{equation}\label{uk_est_fin}
U_k(n) \leq C_4 \cdot  n^{(k+2)\delta_{av}+1}\sum_{1 \leq j \leq k+1} n^{(j+1)(1-\beta)} \leq C_5 \cdot n^{(k+2)(1-\beta+\delta_{av}) + 1}
\end{equation}
for some constants~\(C_4,C_5 > 0.\)

Following an analogous analysis, we get that~\(V_k(n)\) satisfies~(\ref{uk_est_fin}) as well and so
\begin{equation}
J_k(n) \leq C_6 \cdot n^{(k+2)(1-\beta+\delta_{av}) + 1} \label{jk_est}
\end{equation}
for some constant~\(C_6 > 0.\) By Markov inequality~\(\mathbb{P}\left(\sum_{u} d_1(u) d_{k+1}(u) \geq 2J_k(n)\right) \leq \frac{1}{2}\) and so by~(\ref{jk_est}), the sum~\(\sum_{u} d_1(u) d_{k+1}(u)  \leq 2J_k(n) \leq  2C_6 \cdot n^{(k+2)(1-\beta+\delta_{av}) +1}\) with probability at least~\(\frac{1}{2}.\) Together with~(\ref{e_iso}) and~(\ref{m_bound}), we get from~(\ref{nu_low_k}) that
\begin{equation}
\nu_k(G) \geq \frac{m^2}{4\sum_{u} d_1(u) d_{k+1}(u)}  \geq C_7\cdot \frac{n^{4-2\beta+2\delta_{low}}}{n^{(k+2)(1-\beta + \delta_{av}) + 1}} =  C_7\cdot n^{\theta_{low}} \label{nu_low}
\end{equation}
with probability at least~\(\frac{1}{2} - e^{-Cn^{1-\beta+\delta_{low}}} -\exp\left(-\frac{n^{2-\beta + \delta_{low}}}{32}\right) \geq \frac{1}{4}\) for all~\(n\) large, where~\(\theta_{low}\) is as in~(\ref{thet_def}) and~\(C_7 > 0\) is a constant.

Summarizing,
\begin{equation}\label{boot}
\mathbb{P}\left(\nu_k(G) \geq C_7 n^{\theta_{low}}\right) \geq \frac{1}{4} \text{ and so }\mathbb{E}\nu_k(G) \geq \frac{C_7}{4} \cdot n^{\theta_{low}}
\end{equation}
for all~\(n\) large and we use~(\ref{boot}) as a bootstrap to obtain the lower bound in~(\ref{bound_two}). The key ingredient is the following variance estimate:
\begin{equation}\label{var_est_eq}
var\left(\nu_k(G)\right) \leq 4\mathbb{E} \nu_k(G).
\end{equation}
To prove~(\ref{var_est_eq}), we use Theorem~\ref{thm_gen} and let~\(\{Z_e\}_{e \in K_n}\) be independent uniformly distributed random variables in the interval~\([0,1].\) Set the weight of edge~\(e\) as~\(w(e) = \ind(Z_e < p(e)),\) where~\(p(e)\)  is the edge probability as defined in~(\ref{xe_def}). The maximum weight of any edge is at most~\(1\) and so the second moment condition in Theorem~\ref{thm_gen} is satisfied with~\(\mu_2 = 1.\) Consequently, from the first inequality in~(\ref{mk_h_lat}) we get that the variance of~\(\nu_k(G)\) is bounded above by the product of~\(4\mu_2 \leq 4\) and the expected number of edges in a maximum~\(k-\)strong matching of~\(G.\) This proves~(\ref{var_est_eq}).

From~(\ref{var_est_eq}) and Chebychev's inequality we then get for~\(\epsilon > 0\) that
\[\mathbb{P}\left(|\nu_k(G) -\mathbb{E}\nu_k(G)| \geq \epsilon \mathbb{E}\nu_k(G)\right) \leq \frac{var(\nu_k(G))}{(\epsilon\mathbb{E}\nu_k(G))^2} \leq \frac{4}{\epsilon^2 \mathbb{E}\nu_k(G)}\]
and setting~\(\epsilon = \frac{1}{2}\) we get that~\(\mathbb{P}\left(\nu_k(G) \geq \frac{\mathbb{E}\nu_k(G)}{2}\right) \geq 1-\frac{16}{\mathbb{E}\nu_k(G)}.\) Plugging the bound for~\(\mathbb{E}\nu_k(G)\) from~(\ref{boot}), we get the lower bound in~(\ref{bound_two}).~\(\qed\)

\emph{Proof of the upper bound in~(\ref{bound_two})}: The outline is as follows. We first show that\\\(\min_{v} d_k(v) \geq C\cdot n^{k(1-\beta+\delta_{low})}\) for some constant~\(C > 0\) with high probability and then plug this into~(\ref{nu_g_up}) in Proposition~\ref{prop_one} to get the upper bound in~(\ref{bound_two}).

We use an exploration technique analogous to~\cite{boll}. Let~\({\cal S}_0 = \{1\}\) and for~\(i \geq 1\) let~\({\cal S}_i\) be the set of vertices that are at a distance~\(i\) from the vertex~\(1.\) Given~\({\cal S}_0,{\cal S}_1,\ldots,{\cal S}_{i-1} = \{v_1,\ldots,v_L\},\) we would like to estimate the size of~\({\cal S}_i.\)
Let~\({\cal T}_0 = \bigcup_{j=0}^{i-1} {\cal S}_j\) and define the sets~\({\cal T}_l, 1 \leq l \leq L\) iteratively as follows.
Let~\({\cal N}_1\) be the set of neighbours of~\(v_1\) in~\(\{1,2,\ldots,n\} \setminus {\cal T}_0\) and set~\({\cal T}_1 = {\cal T}_0 \cup {\cal N}_1.\)
Iteratively, let~\({\cal N}_j\) be the set of neighbours of~\(v_j\) in the set~\(\{1,2,\ldots,n\} \setminus {\cal T}_{j-1}\) and set~\({\cal T}_j = {\cal T}_{j-1} \cup {\cal N}_j.\)

Define
\begin{equation}\label{p_low_def}
p_{low} := \frac{\gamma_1}{n^{\beta-\delta_{low}}} \text{ and } p_{up} := \frac{\gamma_2}{n^{\beta-\delta_{up}}}
\end{equation}
where the constants~\(\gamma_i,i=1,2\) and~\(\delta_{low},\delta_{up}\) are as in~(\ref{h_cond}). Each edge probability lies between~\(p_{low}\) and~\(p_{up}\) and so if~\(E_j\) denotes the event that
\begin{equation}\label{nj_est}
2p_{up}(n-\#{\cal T}_{j-1}) \geq \#{\cal N}_j \geq p_{low}\cdot \frac{(n-\#{\cal T}_{j-1})}{2},
\end{equation}
then using the standard Binomial deviation estimate~(\ref{conc_est_f}) in Lemma~\ref{lem_one} of Appendix with~\(\epsilon= \frac{1}{2},\) we get that~\(\mathbb{P}(E_j \mid {\cal T}_{j-1}) \geq 1-e^{-C(n-\#{\cal T}_{j-1})p_{low}}\) for some absolute constant~\(C > 0\) not depending on the choice of~\({\cal T}_{j-1}.\)
Thus
\begin{equation}\label{free_hand}
\mathbb{P}\left(\bigcap_{l=1}^{j} E_l \mid {\cal T}_{j-1}\right) \geq \left(1-e^{-C(n-\#{\cal T}_{j-1})p_{low}}\right) \cdot \ind\left(\bigcap_{l=1}^{j-1} E_l\right).
\end{equation}

If~\(\bigcap_{l=1}^{j-1}E_l\) occurs, then~\({\cal T}_{j-1}\) has size at most~\(\#{\cal T}_0 + 2np_{up} (j-1) \leq \#{\cal T}_0 +2np_{up}L\) and so from~(\ref{free_hand}) we get that
\[\mathbb{P}\left(\bigcap_{l=1}^{j} E_l \mid {\cal T}_{j-1}\right) \geq \left(1-e^{-C(n-\#{\cal T}_{0} - 2np_{up}L)p_{low}}\right) \cdot \ind\left(\bigcap_{l=1}^{j-1} E_l\right)\]
and taking expectations and setting~\(j=L,\) we get that
\begin{equation}\label{cruc_est2}
\mathbb{P}\left(\bigcap_{l=1}^{L} E_l \mid {\cal T}_0\right) \geq \left(1-\exp\left(-\frac{Cp_{low}}{4}(n-\#{\cal T}_0 - 2np_{up}L)\right) \right) \cdot \mathbb{P}\left(\bigcap_{l=1}^{L-1} E_l \mid {\cal T}_0\right).
\end{equation}


If~\(\bigcap_{l=1}^{L} E_l \) occurs, then from~(\ref{nj_est}) we see that the number of vertices~\({\cal S}_i\) at a distance~\(i\) from the vertex~\(1\) satisfies
\begin{equation}\label{cruc_est}
2np_{up} L \geq \#{\cal S}_i \geq \sum_{j=0}^{L-1} \frac{p_{low}}{2} \cdot (n-\#{\cal T}_j) \geq \frac{Lp_{low}}{2} \cdot (n-\#{\cal T}_0-2np_{up}L).
\end{equation}
Now define the event~\(J_w:= \left\{\#{\cal S}_{w-1} \cdot \frac{np_{low}}{4} \leq \#{\cal S}_w \leq \#{\cal S}_{w-1} \cdot 2np_{up}\right\}.\) If the event\\\(\bigcap_{1 \leq w \leq i-1}J_w\) occurs, then by iteration we get for~\(1 \leq w \leq i-1\) that
\begin{equation}\label{si_est}
\left(\frac{np_{low}}{4}\right)^{w} \leq {\cal S}_w \leq (2np_{up})^{w}
\end{equation}
and so~\(\#{\cal T}_0 = \sum_{w=1}^{i-1} \#{\cal S}_w \leq \sum_{w=1}^{i-1} (2np_{up})^{w-1} \leq i(2np_{up})^{i-1}.\)
Therefore using the fact that~\(\#{\cal S}_{i-1} = L \leq \#\left(\bigcup_{w=0}^{i-1} {\cal S}_w\right) = \#{\cal T}_0,\)  we get
\begin{equation}\label{tau_est}
\#{\cal T}_0 + 2np_{up}L \leq (2np_{up}+1) \#{\cal T}_0 \leq i (4np_{up})^{i} \leq  k_1 (4np_{up})^{k_1}
\end{equation}
for all~\(n\) large, since~\(i \leq k_1.\) Using the fact that~\((np_{up})^{k_1} \leq C \cdot n^{k_1\left(1-\beta+\delta_{up}\right)}\) for some constant~\(C > 0\) by~(\ref{p_low_def}) and~\(k_1(1-\beta+\delta_{up}) < 1\) strictly (see~(\ref{thet_def})), we get the final term in~(\ref{tau_est})
is bounded above by~\(\frac{n}{8}\) for all~\(n\) large  and so from~(\ref{free_hand}) we see that~\(\mathbb{P}\left(\bigcap_{l=1}^{L} E_l \mid {\cal T}_0\right) \ind \left(\bigcap_{w=1}^{i-1} J_w\right)\) is bounded below by
\begin{equation}
\left(1-\exp\left(-\frac{C}{4} \cdot \frac{np_{low}}{8}\right) \right) \cdot \mathbb{P}\left(\bigcap_{l=1}^{L-1} E_l \mid {\cal T}_0\right) \ind \left(\bigcap_{w=1}^{i-1} J_w\right). \label{iter_j}
\end{equation}
Using the bound~(\ref{iter_j}) iteratively, we then get that
\begin{equation}
\mathbb{P}\left(\bigcap_{l=1}^{L} E_l \mid {\cal T}_0\right) \ind \left(\bigcap_{w=1}^{i-1} J_w\right) \geq \left(1-\exp\left(-\frac{C}{4} \cdot \frac{np_{low}}{8}\right) \right)^{L} \ind \left(\bigcap_{w=1}^{i-1} J_w\right) \label{iter_j2}
\end{equation}
and because~\(L \leq n\) we also have that~\(\left(1-x\right)^{L} \geq 1-Lx\geq 1-nx\) with~\(x =\exp\left(-\frac{C}{4} \cdot \frac{np_{low}}{8}\right) .\)
Thus
\begin{equation}
\mathbb{P}\left(\bigcap_{l=1}^{L} E_l \mid {\cal T}_0\right) \ind \left(\bigcap_{w=1}^{i-1} J_w\right) \geq (1-nx) \ind \left(\bigcap_{w=1}^{i-1} J_w\right). \label{iter_j3}
\end{equation}

Now if~\(\bigcap_{l=1}^{L} E_l \bigcap \bigcap_{w=1}^{i-1}J_w\) occurs, then from~(\ref{tau_est}) and~(\ref{cruc_est}) we get~\[2np_{up} \cdot \#{\cal S}_{i-1} = 2np_{up} L \geq \#{\cal S}_i \] and \[\#{\cal S}_{i} \geq \#{\cal S}_{i-1} \cdot \frac{p_{low}}{2} \cdot (n-\#{\cal T}_0-2np_{up}L) \geq \frac{np_{low}}{4} \cdot \#{\cal S}_{i-1}.\]
In other words, the event~\(J_i\) occurs. From the estimate~(\ref{iter_j3}) we therefore get that
\[\mathbb{P}\left(\bigcap_{w=1}^{i} J_w\right) \geq \left(1-n\cdot e^{-2Dnp_{low}}\right)\cdot \mathbb{P}\left(\bigcap_{w=1}^{i-1} J_w\right) \]
for some constant~\(D >0\) not dependent on the choice of~\(i.\) Proceeding iteratively, we therefore get
\begin{equation}\label{ju_est}
\mathbb{P}\left(\bigcap_{w=1}^{k_1} J_w\right) \geq \left(1-n\cdot e^{-2Dnp_{low}}\right)^{k_1} \geq 1-k_1ne^{-2Dnp_{low}} \geq 1-e^{-D_1n^{1-\beta+\delta_{low}}}
\end{equation}
for all~\(n\) large and some constant~\(D_1 > 0.\)


If~\(\bigcap_{w=1}^{k_1} J_w\) occurs, then from~(\ref{si_est}) we get that~\(d_{k_1}(1) \geq \sum_{w=1}^{k_1}\left(\frac{np_{low}}{4}\right)^{w} \geq \left(\frac{np_{low}}{4}\right)^{k_1}.\)
and so from~(\ref{ju_est}) and the union bound, we get that
\begin{equation}\label{dk_min}
\min_{v} d_{k_1}(v) \geq \left(\frac{np_{low}}{4}\right)^{k_1} \geq C \cdot n^{k_1(1-\beta+\delta_{low})}
\end{equation}
for some constant~\(C > 0\) with probability at least~\(1-n\cdot e^{-D_1n^{1-\beta+\delta_{low}}} \geq 1-e^{-D_2n^{1-\beta+\delta_{low}}}\) for all~\(n\) large and some constant~\(D_2 > 0.\) Plugging~(\ref{dk_min}) into the upper bound~(\ref{nu_g_low}) in Lemma~\ref{prop_one} of Appendix, we get the upper bound in~(\ref{bound_two}).~\(\qed\)

\setcounter{equation}{0}
\renewcommand\theequation{A.\arabic{equation}}
\section*{Appendix}
The following result obtains bounds on the strong matching number of any graph~\(H\) in terms of local vertex neighbourhood sizes.
\begin{lemma}\label{prop_one} Let~\(H =(V,E)\) be any graph containing no isolated edge and for~\(j \geq 1,\) let~\(d_{j}(u)\) denote the number of vertices that are at a distance at most~\(j\) from the vertex~\(u \in V\) and let~\(\Delta  = \Delta(H) := \max_{u \in V} d_1(u)\) be the maximum degree of a vertex in~\(H,\) respectively. If~\(H\) contains~\(n\) vertices and~\(m\) edges, then for~\(k \geq 0\) we have
\begin{equation}\label{nu_g_low}
\nu_k(H) \geq \frac{m^2}{4} \left(\sum_{u \in V} d_{1}(u)(d_{k+1}(u)-1)\right)^{-1} \geq \frac{m}{8\Delta^{k+1}},
\end{equation}
and for~\(k \geq 3\) we have that
\begin{equation}\label{nu_g_up}
\nu_k(H) \leq \frac{n}{\min_{u}d_{k_1}(u)},
\end{equation}
where~\(k_1= \frac{k-1}{2}\) if~\(k_1\) is odd and~\(k_1 = \frac{k-2}{2}\) otherwise.
\end{lemma}
The first expression of~(\ref{nu_g_low}), which we call the \emph{average degree bound}, provides a lower bound in terms of average neighbourhood size of vertices and is particularly useful in the context of inhomogenous random graphs, (see Section~\ref{inhom_strong}). From the second bound in~(\ref{nu_g_low}), which we call the \emph{maximum degree bound}, we see that graphs with bounded maximum degree have a~\(k-\)strong matching number that grows at least linearly with the number of edges. The upper bound in~(\ref{nu_g_up}) implies that the strong matching number is low for graphs with large minimum degree, since in such graphs, the edges are connected to each other through paths of small lengths.

We prove the upper and lower bounds for~\(\nu_k(H)\) in Lemma~\ref{prop_one} in that order.\\
\emph{Proof of~(\ref{nu_g_up}) in Lemma~\ref{prop_one}}: Let~\({\cal F} = \{(u_1,v_{1}),\ldots,(u_t,v_{t})\}, t=\nu_k(H)\) be a maximum~\(k-\)strong matching in~\(H.\) Let~\({\cal N}_j(u)\) be the set of all vertices at a distance at most~\(j\) from~\(u\) in the graph~\(H\) and so letting~\(k_1\) be as defined in the statement of the Theorem, we have that~\(\bigcup_{i=1}^{t} {\cal N}_{k_1}(u_i) \subseteq V(H).\) For any two vertices~\(u_i\) and~\(u_j\) belonging to distinct edges of the matching~\({\cal F},\) we must have that~\({\cal N}_{k_1}(u_i) \bigcap {\cal N}_{k_1}(u_j) = \emptyset;\) because otherwise, there would be a path of length at most~\(2k_1 \leq k-1\) connecting~\(u_i\) and~\(u_j.\)
Therefore~\(t \cdot \min_{u} \#{\cal N}_{k_1}(u) \leq \# V(H) = n\) and this obtains~(\ref{nu_g_up}).~\(\qed\)

\begin{figure}[tbp]
\centering
\includegraphics[width=3in, trim= 20 400 50 130, clip=true]{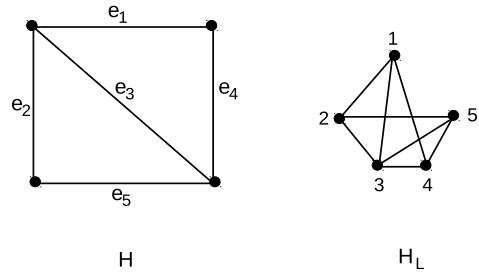}
\caption{A graph~\(H\) (left) and its corresponding line graph~\(H_L\) (right) where vertex~\(i\) in~\(H_L\) corresponds to the edge~\(e_i\) in~\(H.\) The edges~\(e_1\) and~\(e_5\) form a matching in~\(H\) and the corresponding vertices~\(1\) and~\(5\) form a stable set in~\(H_L.\) }
\label{fig_squares}
\end{figure}

\emph{Proof of~(\ref{nu_g_low}) in Lemma~\ref{prop_one}}: Let~\(H_L\) be the line graph (pp. 71,~\cite{har}) of~\(H\) obtained by assigning a vertex~\(v_e\) for each edge~\(e \in H\) and connecting~\(v_e\) and~\(v_f\) by an edge in~\(H_L\) if and only if the edges~\(e\) and~\(f\) share a common endvertex  in~\(H.\) If~\({\cal I} = \{v_{f_1},\ldots,v_{f_t}\}\) is a stable vertex set in~\(H_L,\) i.e. no two vertices in~\({\cal I}\) are adjacent to each other in~\(H_L,\) then the corresponding edges~\(f_1,\ldots,f_t\) form a matching in~\(H.\) In Figure~\ref{fig_squares} we have illustrated the above reasoning with an example, where the edges~\(e_1\) and~\(e_5\) form a matching. The corresponding vertices~\(1\) and~\(5\) in~\(H_L\) form a stable set.

Let~\(H^{k+1}_L\) be the~\((k+1)^{th}\) power of the line graph~\(H_L\) obtained as follows: The vertex set of~\(H^{k+1}_L\) is the same as that of~\(H_L.\) If the distance between the vertices~\(v_e\) and~\(v_f\) is at most~\(k+1\) in~\(H_L,\) then~\(v_e\) and~\(v_f\) are adjacent in~\(H^{k+1}_L.\) By construction, if~\(\{v_{e_1},\ldots,v_{e_r}\}\) is a stable set in~\(H^{k+1}_L,\) then the set of edges~\(\{e_1,\ldots,e_r\}\) form a~\(k-\)strong matching in the graph~\(H.\)

The number of vertices in~\(H^{k+1}_L\) equals~\(m,\) the number of edges in~\(H.\) Moreover if~\(d_{av}\) denotes the average degree of a vertex in~\(H^{k+1}_L,\) then using the fact that no edge of~\(H\) is isolated, we get that~\(d_{av} \geq 1.\) From Theorem~\(3.2.1,\) pp.~\(29,\)~\cite{alon}, we therefore get that there exists a stable set in~\(H_L^{k+1}\) of size at least~\(\frac{m}{2d_{av}}.\) To estimate~\(d_{av},\) let~\({\cal N}_{j}(u)\) denote the set of vertices at a distance at most~\(j\) from the vertex~\(u \in H.\) The degree of the vertex~\(v_e\) in the graph~\(H^{k+1}_L\) corresponding to the edge~\(e = (u,v) \in H\) is at most~\(\#\left(\left({\cal N}_{k+1}(u) \setminus \{v\}\right)  \bigcup \left({\cal N}_{k+1}(v) \setminus \{u\}\right)  \right)\) and so
\begin{eqnarray}
d_{av} &\leq& \frac{1}{m} \sum_{e = (u,v) \in H} \#\left(\left({\cal N}_{k+1}(u) \setminus \{v\}\right)  \bigcup \left({\cal N}_{k+1}(v) \setminus \{u\}\right)  \right) \nonumber\\
&\leq& \frac{1}{m} \sum_{e = (u,v) \in H} \#\left({\cal N}_{k+1}(u) \setminus \{v\}\right) + \#\left({\cal N}_{k+1}(v) \setminus \{u\}\right) \nonumber\\
&=& \frac{2}{m} \sum_{u \in V} \sum_{v \in {\cal N}_1(u) \setminus \{u\}}  \#\left({\cal N}_{k+1}(u) \setminus \{v\}\right) \nonumber\\
&=& \frac{2}{m} \sum_{u \in V} d_{1}(u)(d_{k+1}(u)-1). \label{d_av_temp}
\end{eqnarray}
This obtains the first lower bound in~(\ref{nu_g_low}).

Finally, if~\(\Delta\) is the maximum degree of a vertex in~\(H,\) then~\(d_{k+1}(u) \leq \Delta^{k+1}\) and so from~(\ref{d_av_temp}) we get that~\(d_{av} \leq \frac{2\Delta^{k+1}}{m} \sum_{u \in V} d_1(u)  = 4\Delta^{k+1}.\) This proves the second lower bound in~(\ref{nu_g_low}).~\(\qed\)

Also, throughout we use the following standard deviation estimate.
\begin{lemma}\label{lem_one}
Let~\(\{X_j\}_{1 \leq j \leq L}\) be independent Bernoulli random variables with~\(\mathbb{P}(X_j = 1) = 1-\mathbb{P}(X_j = 0) > 0.\) If~\(T_L = \sum_{j=1}^{L} X_j,\theta_L = \mathbb{E}T_L\) and~\(0 < \epsilon \leq \frac{1}{2},\) then
\begin{equation}\label{conc_est_f}
\mathbb{P}\left(\left|T_L - \theta_L\right| \geq \theta_L \epsilon \right) \leq 2\exp\left(-\frac{\epsilon^2}{4}\theta_L\right).
\end{equation}
\end{lemma}
For a proof of Lemma~\ref{lem_one}, we refer to Corollary A.1.14, pp. 312 of~\cite{alon}.

\ACKNO{I thank Professors Rahul Roy, C. R. Subramanian and the referees for crucial comments that led to an improvement of the paper. I also thank IMSc for my fellowships.}


\end{document}